\documentstyle{article}
\newcommand{\R}{{\bf R}}
\newcommand{\length}{{\rm length}} 
\newcommand{\size}{{\rm size}}

\newcommand{\Z}{{\bf Z}}

\newcommand{\CP}{{\bf CP}}

\newcommand{\eps}{{\varepsilon}}

\newcommand{\cR}{{\cal R}}

\newcommand{\Ham}{{\rm Ham }}
\newcommand{\Symp}{{\rm Symp}}

\def\one{\hbox{1\hskip-2.4pt l}}

\newcommand{\niN}{{\noindent}}

\newcommand{\QED}{\hfill$\Box$\medskip}

\title{Symplectic aspects of the first eigenvalue}
\author{Leonid Polterovich\thanks{Supported by  
the United States - Israel Binational Science Foundation 
grant 94-00302
} 
\\ Tel-Aviv University
\\ (polterov \@@ math.tau.ac.il) }

\begin{document}
\maketitle

\bigskip

There are two themes in the present paper. The first one
is spelled out in the title, and is inspired by an attempt to find an
analogue of Hersch-Yang-Yau estimate for $\lambda_1$ of surfaces
in symplectic category. In particular we prove
that every split symplectic manifold $T^4 \times M$ admits
a compatible Riemannian metric whose first eigenvalue is arbitrary
large. On the other hand for K\"ahler metrics compatible with a given
integral symplectic form an upper bound for $\lambda_1$ does exist.
The second theme is the study of 
Hamiltonian symplectic fibrations over $S^2$. 
We construct a numerical invariant called {\it the size of 
a fibration} which arises as the solution of certain variational
problems closely related to Hofer's
geometry, K-area and coupling. In some examples it can be computed 
with the use of Gromov-Witten invariants. The link between these two
themes is given by an observation that the first eigenvalue of a
Riemannian metric compatible with a symplectic fibration
admits a universal upper bound in terms of the size.

\bigskip

{\bf 1. Introduction and main results}

\medskip

{\it 1.1 Motivation and an overview}

 Let $M$ be a smooth closed manifold
endowed with some "geometric" structure. In many cases this structure
determines in a canonical way a class $\cR$ of Riemannian metrics
on $M$. Traditionally, one asks the following

\medskip

{\bf Question.} Does there exist {\it a universal}
upper bound $\lambda_1(M,g) \leq C$ for all $g \in \cR$?

\medskip

Here and below we write $\lambda_1(M,g)$ for the first
positive eigenvalue of the Laplacian of $(M,g)$ acting
on functions. Let us give several examples.

\medskip

{\it 1.1.A Topology}(no additional structure on $M$).
Here $\cR$ is the class of all metrics on $M$ whose
sectional curvatures at each point belong to $[-1;1]$.
In this case existence of such a universal bound leads
to some non-trivial topological restrictions on $M$
(Wu [W], 1995).

\medskip

{\it 1.1.B Surfaces endowed with an area form}.
Let $\cR$ be the class of all metrics on an orientable
closed surface $M$ with a given area form. Then for all
$g \in \cR$ holds
$$\lambda_1(M,g) \leq {{{8\pi({\rm genus}(M) + 1)}\over
{{\rm Area}(M)}}}.$$
This was proved by Hersch [He], 1970 for $M=S^2$ and
by Yang and Yau [Y-Y], 1980 for surfaces of higher genus.

The beauty of this result inspired many attempts to find
its multi-dimensional analogue (see Gromov [G1], 1993
for a stimulating discussion). The difficulty here is {\it to
recognize "to which category" belongs this estimate}. Indeed,
in dimension 2 one cannot distinguish between a volume form
and a symplectic structure, and between Riemannian and K\"ahler
geometry. Some important
developments in this direction are 
reflected in the next two examples.

\medskip

{\it 1.1.C Higher-dimensional manifolds endowed
with a volume form.}
Let $\cR$ be the class of all Riemannian metrics
on $M$ with a given volume form. In 1993 Colbois
and Dodziuk [C-D] proved that $\sup_{g\in \cR}\lambda_1(M,g) = \infty$
for every $M$ of dimension $\geq 3$.

\medskip

{\it 1.1.D Complex manifolds with a distinguished cohomology class.}
Let $\cR$ be the class of all K\"ahler metrics on a complex
manifold $M$ whose K\"ahler form represents a given 
cohomology class. In various interesting cases a universal upper bound 
for $\lambda_1(M,g),\; g \in \cR$ exists and is related to
algebraic-geometric properties of $M$. This was discovered
by Li and Yau [L-Y] in 1982, and by Bourguignon, Li and Yau
[B-L-Y] in 1994 (see 1.2, 1.5 and 4.1 below for more detailes).

\medskip

In the present paper we study non-existence/existence
of universal upper bounds for $\lambda_1$ in the symplectic
context, which is a natural generalization of an area
form on a oriented surface. We fix a closed symplectic
manifold $(M,\Omega)$ and consider certain classes 
$\cR$ of {\it compatible} Riemannian metrics, that is metrics
$g$ of the form
$$g(\xi,\eta)=
\Omega(\xi,J\eta),\; \xi,\eta \in TM,$$
where $J$ is an almost complex structure on $M$.

At the beginning 
we focus on two cases: $\cR$ consists of {\it all} compatible
metrics, and $\cR$ consists of {\it K\"ahler} compatible 
metrics (that is metrics associated to integrable $J$'s).
In a number of examples we establish "flexibility"
(that is $\sup \lambda_1 = \infty$) in the first case (see 1.2.A below)
vs. "rigidity" ($\sup \lambda_1 < \infty$) in the second one (see 
1.2.B below). Let us mention that our point of view on the K\"ahler
case is different from one in 1.1.D: we fix $\Omega$ and vary $J$,
while in 1.1.D the complex structure $J$ is fixed and the compatible
 symplectic
form $\Omega$ varies in a given cohomology class.

However, this is not the end of the story. It turns out that
there exists an intermediate (with respect to "degree of integarbility"
of $J$'s) category where in many situations an upper bound
for the first eigenvalue does exist, namely 
{\it symplectic fibrations over $S^2$} endowed with a class of
Riemannian metrics which satisfy a natural compatibility condition
(see
1.4
and 1.5
below). Our upper bound for $\lambda_1$ of a symplectic fibration turns
out to be related to such objects of Symplectic Topology as
Hofer's metric, K-area and coupling (see 1.9 below). In this aspect
the present research continues the study of "hard" numerical
invariants of symplectic fibrations initiated in [P1-P3].

\bigskip

{\it 1.2 $\lambda_1$ for a symplectic manifold}

Let $(M,\Omega)$ be a symplectic manifold.
Set $\lambda_1(M,\Omega) = \sup_g \lambda_1(M,g)$,
where $g$ runs over all compatible metrics. As we have seen in
1.1.B the quantity $\lambda_1 (M,\Omega)$ is finite
when $M$ is a closed surface. Our first result states that
at least {\it after suitable stabilization} this property
disappears and hence the situation becomes flexible as in 1.1.C.

More precisely, consider the 4-torus $T^4 = \R^4/\Z^4$ endowed
with the standard symplectic form $\sigma = dp_1 \wedge dq_1 +
dp_2 \wedge dq_2$, where $(p_1,q_1,p_2,q_2)$ are coordinates
on $\R^4$.

\proclaim Theorem 1.2.A. Let $(M,\Omega)$ be a closed symplectic
manifold. Then 
$$\lambda_1(T^4 \times M, \sigma \oplus \Omega) = \infty.$$

Notice that this result remains true when $M$ is a point.
The proof which is given in section 2 below is based
on the construction of an auxiliary hypoelliptic operator
in the spirit of [B-B]. It would be interesting to 
understand to which extent the stabilization is important.

The situation changes drastically when one considers
a smaller class of K\"ahler metrics whose K\"ahler
class equals to $\Omega$. 

\proclaim Theorem 1.2.B. Let $(M,\Omega)$ be a closed
symplectic manifold of real dimension $2n$
such that the cohomology class
of $\Omega$ is integral: $[\Omega] \in H^2(M,\Z)$.
Then for every compatible K\"ahler metric $g$ on $M$
holds
$$\lambda_1(M,g) \leq {\rm constant} (n+2 -
{{(c_1(TM)\cup [\Omega]^{n-1},[M])}\over
 {([\Omega]^n,[M])}}),$$
where the constant depends only on $n$.

\medskip
\niN
The proof which is a simple combination of two known
results goes as follows. Given a K\"ahler projective
manifold, one can estimate $\lambda_1$ in terms of its
degree [B-L-Y]. On the other hand, there exists
an upper bound for the projective degree of an integral K\"ahler
manifold in purely cohomological terms. This bound which
can be considered as an effective version of the classical
Kodaira embedding theorem was established 
by Demailly [De]. The details of this argument can be found
in section 4.

{\bf Remark 1.2.C} Clearly, 1.2.B automatically implies an upper bound
for symplectic manifolds whose symplectic form
represents a {\it rational} cohomology class,
that is $[\Omega] \in {1\over N}H^2(M,\Z)$, where $N$ 
is an integer. However this upper bound goes to infinity with $N$.
Therefore the approximation argument does not work, and
the situation in the case when $[\Omega]$ is 
{\it irrational} remains completely unclear even for
simplest manifolds (for instance for the 4-dimensional torus).

\bigskip

{\it 1.3 Symplectic fibrations over $S^2$}

Let $(M, \Omega)$ be a closed symplectic manifold.
A fibration $p:P \to S^2$ with fiber $M$ over an oriented 2-sphere
is called {\it symplectic}
if each fiber $p^{-1}(x)$ is endowed with a symplectic form
$\Omega_x$ which depends smoothly on the base point $x$,
and is symplectomorphic to $\Omega$. Such fibrations
which were introduced
in [G-L-S] appear in various problems of Symplectic Topology,
in particular in the study of loops on the group of symplectomorphisms.
We are going now to explain this link in more details.

Let $\{f_t \},\; t \in S^1, \; f_0 = f_1 = \one$ be a loop
of symplectic diffeomorphisms of $(M, \Omega)$
which represents an element
$\gamma \in \pi_1(\Symp (M,\Omega),\one )$.
Let $D_+$ and
$D_-$ be two copies of the disc $D^2$ bounded by $S^1$.
Consider a map $\Psi : M \times S^1 \to M \times S^1$ 
given by $(z,t) \to (f_t z,t)$. Define now a new manifold
$$P({\gamma}) = (M \times D_-) \cup_{\Psi} (M \times D_+).$$
It is clear that $P_{\gamma}$ has the canonical structure
of a symplectic fibration over $S^2$ which depends only
on the element $\gamma$. In what follows we assume that the
base $S^2$ is oriented, and the orientation comes from $D_+$.
Note that the choice of an orientation on $S^2$ is equivalent
to the choice of an orientation on $P$.

Clearly, this construction can be reversed. Namely,
given a symplectic fibration $p:P \to S^2$ over
an oriented 2-sphere with the
fiber $(M, \Omega)$, one can reconstruct the corresponding
element $\gamma \in \pi_1(\Symp(M,\Omega))$ up to conjugation.
For that purpose, fix a point, say $x$, on $S^2$ and
consider  a loop $\alpha(s),\; s\in S^1$ in the
space $LS^2$ of all loops on $S^2$ based at $x$.
Assume that $\alpha(0)$ is the 
constant loop, and that $\alpha$ represents a {\it positive}
generator in $\pi_1(LS^2)  = \pi_2 (S^2)$ (here we use
the choice of an orientation on $S^2$). Choose a symplectic
connection on $P$ (we refer the reader to [G-L-S], [McD-S]
for basic notions of the theory of symplectic connections).
Then the holonomies of this connections over loops $\alpha(s)$
form a loop in the group of symplectomorphisms of the fiber
$P_x$ over $x$, which defines the needed element $\gamma$.

Throughout this paper, we restrict our attention
to the subgroup $$\Ham (M,\Omega) \subset \Symp (M,\Omega)$$
of all Hamiltonian diffeomorphisms of $(M,\Omega)$.
In particular we consider only Hamiltonian loops and
corresponding symplectic fibrations, which we call
{\it Hamiltonian symplectic fibrations}. Hamiltonian
fibrations can be characterized by the fact that 
there exists a cohomology class on the total space
whose restriction to fibers coincides with the class
of the symplectic form (see [Se]) . Notice that the difference between
"symplectic" and "Hamiltonian" disappears when the manifold $M$
is simply connected.

\bigskip

{\it 1.4 $\lambda_1$ of a symplectic fibration over $S^2$}

Let $p :P \to S^2$ be a Hamiltonian symplectic fibration.
A pair $(\omega, j)$ consisting of a symplectic form
$\omega$ on $P$ and an almost complex structure $j$
on $P$ is called a {\it compatible quasi-K\"ahler 
structure} if the following holds:

(i) The form
$\omega (\xi,j\eta)\;,\xi,\eta \in TP$ is a
Riemannian metric on $P$;

(ii) The restriction of $\omega$ to fibers coincides
with the symplectic structure on fibers;

(iii) The projection $p$ is $(j,i)$-holomorphic
for some complex structure $i$ on $S^2$ which respects
the orientation.

\medskip

We say that a Riemannian metric $g$ on $P$ is {\it compatible}
if it comes from a compatible quasi-K\"ahler structure.

A natural way to construct compatible  quasi-K\"ahler
structures on $P$ is as follows. Fix a field $J_x,\; x \in S^2$
of almost complex structures on the fibers of $P$ which
are compatible with symplectic forms on the fibers. Pick
up a complex structure $i$ on $S^2$.
Let $\omega$ be a
symplectic form on $P$ which extends symplectic structures
on the fibers and defines the positive orientation on $P$.
The field of $\omega$-orthogonal subspaces to the fibers
defines a connection on $P$, that is a splitting
$T_{(x,z)}P = T_zP_x \oplus T_xS^2$ where $x \in S^2$,
$P_x$ is a fiber over $x$ and $z \in P_x$. Consider
an almost complex structure $j = J_x \oplus i$ on $P$.
Clearly, $(\omega, j)$ satisfies (i) - (iii) above.

Now we are ready to define a spectral invariant of
a symplectic fibration. Namely, set
$$\lambda_1(P) = \sup_g \lambda_1(P,g),$$
where the supremum is taken over all compatible
Riemannian metrics on $P$. 

It turns out that in some interesting examples
this quantity is finite. 
Consider for instance the case
$(M,\Omega)= (\CP^{n-1},\Omega_{st})$, where $\Omega_{st}$
is the Fubini-Study form normalized in such a way that the
integral over a projective line equals to 1. Let
$\gamma_{q,n} \in \pi_1(\Ham(\CP^{n-1},\Omega_{st}))$
be the element represented by the $S^1$-action
$$f_{q,n}: (z_1:...:z_n) \to
(e^{2\pi it}z_1:...:e^{2\pi it}z_q:z_{q+1}:...:z_n),$$
where $q \in \{1, ..., n-1\}$. Let $P_{q,n}$ be the corresponding
symplectic fibration. 

\proclaim Theorem 1.4.A. The following inequality
holds:
$$\lambda_1(P_{q,n}) \leq {{8\pi n}\over {n-q}}.$$

\medskip

Moreover, this estimate is {\it stable } in the following sense.
Let $(T^4,\sigma)$ be the standard symplectic 4-torus as in
1.2. Take the product
$\tilde P_{q,n} = T^4 \times P_{q,n}$ which we consider
as a symplectic fibration over $S^2$ with the fiber 
$$(T^4 \times \CP^{n-1}, \sigma \oplus \Omega_{st}).$$

\proclaim Theorem 1.4.A'. The following inequality
holds:
$$\lambda_1(\tilde P_{q,n}) \leq {{8\pi n}\over {n-q}}.$$

The proof of both estimates is given in 1.9 below.
It seems to be a difficult problem to compute
this invariant even in simple examples. Consider
for instance the simplest case of the trivial
symplectic fibration with the fiber $(M,\Omega)$,
that is $P_{trivial} = M \times S^2$.
Clearly, in this situation we can consider
{\it split} compatible metrics, and  making the area 
of the base $S^2$ arbitrarily small we get that
$$\lambda_1(P_{trivial}) \geq \lambda_1 (M, \Omega).$$
Applying results of 1.2 we see that $\lambda_1$ is
infinite for trivial fibrations with a stabilized
fiber of the form $T^4 \times M$. Comparing this with 1.4.A'
above we come to the following conclusion:

\niN
{\it invariant $\lambda_1(P)$ distinguishes
the trivial fibration with the fiber $T^4 \times \CP^{n-1}$
from the non-trivial fibrations $\tilde P_{q,n}$.}

\niN 
Without stabilization we are able to prove the same statement
under an additional assumption that 
$n \geq 4$ and $q \leq n-3$. Indeed, 
write $P_{trivial}$ for the trivial symplectic
fibration with the fiber $\CP^{n-1}$.
Arguing as above we see that
$$\lambda_1(P_{trivial}) \geq \lambda_1(\CP^{n-1}, {\rm
Fubini-Studi}) = 4\pi n.$$
On the other hand, applying 1.4.A we get that
$$\lambda_1(P_{q,n}) \leq {8\pi n \over {n-q}} < 4\pi n \leq
\lambda_1(P_{trivial}),$$
and hence our spectral
invariant distinguishes $P_{q,n}$ from the trivial fibration.

\medskip

\niN
{\bf Remark 1.4.B.} It is instructive to compare 1.4.A,A' above
with the results of 1.2. As we have seen in 1.2.B the integrability
of a compatible almost complex structure imposes a global 
restriction on the geometry, namely it prevents the first eigenvalue
from being arbitary large. The almost complex structures considered
in the present section are "integrable in one direction", namely
by definition they admit a holomorphic mapping to $\CP^1$. The
message of theorems 1.4.A,A' above (and of the more general statement
1.5.B below) is that this small amount of integrability is still 
sufficient for 
the global geometric restriction of the same nature. However
theorem 1.2.A shows that it disappears if one considers 
general non-integrable almost complex structures. 

\vfill\eject

{\it 1.5 The Li-Yau estimate revisited.}

The first step towards getting an upper bound
for $\lambda_1(P)$ is the following minor generalization
of a theorem by Li and Yau [L-Y].
Let $(\omega, j)$ be a compatible quasi-K\"ahler structure
on a Hamiltonian symplectic fibration $P$ with a
fiber $(M,\Omega)$, and let $g$
be the corresponding Riemannian metric. 
Set
$$s(\omega) = {{{\rm Volume}(M,\Omega)} \over {{\rm Volume}
(P,\omega)}}.$$
Then $\lambda_1(P,g) \leq 8\pi s(\omega)$. The proof of this
inequality goes exactly as in [L-Y, Theorem 3] (let us 
mention that our normalization conventions are slightly
different from the ones in [L-Y]). Note that in our situation
$j$ is not assumed to be integrable on the fibers, however
this does not play any role.

\medskip
\niN
{\bf Definition 1.5.A} Let $P \to S^2$ be
a Hamiltonian symplectic fibration. Define
its {\it size} as follows:
$$\size (P) = \sup_{\omega} s(\omega),$$
where $\omega$ runs over all symplectic forms
on $P$ which extend the symplectic structures on the fibers and
define the positive orientation.

\medskip

With this language our previous discussion leads to the
next statement.

\proclaim Theorem 1.5.B. The following inequality holds:
$$\lambda_1(P) \leq 8\pi \size(P).$$

Now we face the following purely symplectic question.
{\it Given a symplectic fibration, is its size 
finite ?} It turns out that $\size (P)$ admits three
other equivalent definitions in more familiar terms.  Namely,
one can describe it using
Hofer's geometry, K-area and coupling (see 1.9.A below).
Having in mind this description, we get the finitness of the size
(and in fact compute it explicitely) in a number of examples
using results of [P1].

\medskip
{\it 1.6 Hofer's geometry}

Given a path $\{f_t\}$ in $\Ham (M,\Omega) $ ($t \in [0;1]$),
define its {\it length} as follows. Let $F: M \times S^1 \to \R$ be 
the associated
Hamiltonian function normalized so that the mean value
of $F_t$ over $M$ vanishes for all $t \in S^1$. Then
$\length \{f_t\} = \max |F|$. With this language the distance between
two Hamiltonian diffeomorphisms is just the minimal possible length
of a path between them (see [H-Z], [McD-S] for an introduction
to Hofer's geometry). One can show
that this distance, which was called
{\it coarse} Hofer's distance in [P1], in fact coincides with
the one associated to the "maximum of the
absolute value" norm on the Lie algebra of all normalized
Hamiltonian functions.

In the present paper, we deal mainly not with the Hofer's
length itself, but with {\it its positive and negative parts}
(see [E-P, 4.4]). Namely set
$\length_+ \{f_t\} = \max F$ and
$\length_- \{f_t\} =- \min F$. 

For an element $\gamma \in \pi_1(\Ham)$ define its norm
as 
$$||\gamma|| = \inf \; \length\{f_t\},$$
where the infimum is taken over all loops $\{f_t\}$
with $f_0 = f_1 = \one$
representing $\gamma$. Define also its positive and negative
parts as
$$||\gamma||_+ = \inf \; \length_+\{f_t\},$$ and

$$||\gamma||_- = \inf \; \length_-\{f_t\}.$$ 

Clearly, $||\gamma|| \geq \max(||\gamma||_+,||\gamma||_-)$.
Interestingly enough, in all known examples in fact the
equality take place, however we cannot prove this in general
situation.

Note also that $||\gamma||_- = ||\gamma ^{-1}||_+$
(see 1.7.A below for further discussion about this duality).

\bigskip

{\it 1.7 K-area of symplectic fibrations}

Let $p :P \to S^2$ be a Hamiltonian symplectic fibration
with fiber $(M,\Omega)$. A connection $\nu$ on $P$ is
called symplectic if the parallel transport preserves
symplectic forms on fibers. The curvature $\rho^{\nu}$ of 
$\nu$ is a 2-form on the base which take values in
the Lie algebra of the group of symplectic diffeomorphisms
of a fiber, that is in the space of locally Hamiltonian vector
fields. A symplectic connection is called Hamiltonian if its
curvature takes values in {\it Hamiltonian vector fields}.
In this case we make a natural identification and consider
$\rho^{\nu}(\xi,\eta)$ as a Hamiltonian function on the fiber
with zero mean value. We call such Hamiltonian functions 
{\it normalized}. It is not hard to show that a symplectic
fibration is Hamiltonian in the sense of 1.3 if and only if
it admits a Hamiltonian symplectic connection.

Given an area form, say $\tau$ on $S^2$, one can write
$\rho^{\nu} = L^{\nu}\tau$, where $L^{\nu}$ is a function on $P$.

Fix an area form $\tau$ on $S^2$ of the total area 1
which respects the orientation (recall, that $S^2$ is always
assumed to be oriented). Define the K-area of $P$ as follows ([P1], cf.
[G3]):

$$\chi (P) = \sup_{\nu}{1\over{\max_P |L^{\nu}|}}\;,$$

\niN
where the supremum is taken over all Hamiltonian symplectic connections
$\nu$ on $P$. Notice that the definition in [P1] 
contained
an additional
multiple $1\over {2\pi}$.

 In what follows we shall consider also
{\it positive and negative parts} of the K-area:

$$\chi_+ (P) = \sup_{\nu}{1\over{\max_P L^{\nu}}}\;$$

and

$$\chi_- (P) = \sup_{\nu}{1\over{-\min_P L^{\nu}}}\;.$$

Note that $\chi(P) \leq \min(\chi_-(P),\chi_+(P))$,
and in all known examples in fact the equality holds.
Note also that {\it all three quantities
do not depend on the specific choice of the
area form $\tau$}. 

\niN
{\bf Remark 1.7.A $\;$ (Duality).}
Let $P \to S^2$ be a Hamiltonian symplectic fibration over
an oriented 2-sphere
associated to an element $\gamma \in \pi_1(\Ham(M,\Omega))$.
Change the orientation of the base $S^2$ and denote the
new object by $P^*$. Clearly $P^*$ as a Hamiltonian
symplectic fibration over the oriented 2-sphere
is associated with the element $\gamma^{-1}$.
Note that the following identity holds:
$\chi_-(P) = \chi_+(P^*)$.
 
\medskip

{\it 1.8 Coupling}

Let $p : P \to S^2$ be a symplectic fibration associated to
an element
$\gamma \in \pi_1(\Ham(M,\Omega))$. 
Denote by $a$ the positive generator
of $H^2(S^2,\Z)$. Let
$c \in H^2(P,\R)$ be the unique class whose restriction on fibers
coincides with $[\Omega]$ and whose top power vanishes. Following
[G-L-S] we call $c$ {\it the coupling class} of $P$.

The {\it weak coupling} construction [G-L-S]
prescribes
that for a sufficiently small $\eps > 0$ there exists a 
a smooth family of closed 2-forms
$\{\omega_t\}, \; t \in [0;\eps)$ on $P$
with the following properties:

(i) $\omega_0$ is the lift of an area form on $S^2$
(the total area equals to 1);

(ii)$[\omega_t] = tc + [p^*a]$;

(iii) the restriction of $\omega_t$ to each fiber of
$P$ coincides with a multiple of the symplectic form on the fiber;

(iv) $\omega_t$ is symplectic for $t > 0$.

Following [P1], define $\epsilon(P)$ as the supremum
of all such $\epsilon$ (this invariant measures how strong
the weak coupling can be). 

\medskip

{\it 1.9 On the size of a symplectic fibration}

Now we are ready to state our main result on the
size of a symplectic fibration.

\proclaim Theorem 1.9.A. Let $P \to S^2$ be
a symplectic fibration associated to an element
$\gamma \in \pi_1(\Ham(M,\Omega))$. The following
identity holds:
$$ \size(P) = \epsilon(P) = \chi_+(P) = {1\over{||\gamma||_+}}.$$

\medskip

The proof is given in section 3 below.

\medskip
\niN
{\bf Remark 1.9.B}
We will prove also that
$\chi(P) = {1\over{||\gamma||}}$ (see 3.1.B; note that the
inequalities 
$\chi(P) \geq {1\over{||\gamma||}}$ 
and $\chi(P) \leq \epsilon(P)$ were proved in [P1]).
This identity in a sense reflects a priniciple of the Yang-Mills theory
that the critical values of the Yang-Mills functional on a given 
$G$-bundle over $S^2$ correspond to the lengths of closed geodesics on
the Lie group $G$, see [A-B],[Gra]. Our proof in 3.3 is somewhat 
similar to the argument in [Gra].

\medskip

Clearly, the size of the trivial fibration is infinite.
The general phylosophy of Hofer's geometry suggests
that $||\gamma||_+$ is strictly positive for non-trivial
elements $\gamma$, and hence one can hope that
non-trivial symplectic fibrations has finite size.
In a number of situations this was confirmed in [P1],[P2]
with the use of Gromov-Witten invariants. In particular,
one can use these results in order to prove 1.4.A and 1.4.A'.

\medskip
\niN
{\bf Proof of 1.4.A and 1.4.A':} We use notations of 1.4.
It was shown in [P1] that
$\epsilon(P_{q,n}) = {n\over{n-q}}$. Moreover, the argument
of [P1] remains valid under stabilization, thus
$\epsilon(\tilde P_{q,n}) = {n\over{n-q}}$. 
Thus 1.9.A implies that

$$\size(P_{q,n}) = \size (\tilde P_{q,n}) = {n\over{n-q}}.$$

The desired statements follow immediately from 1.5.B.
\QED

\medskip

The rest of the paper is organized as follows.
In section 2 we prove Theorem 1.2.A, in section
3 we prove Theorem 1.9.A, and in section 4 we prove Theorem 1.2.B.

\bigskip

{\bf 2. Collapsing symplectic manifolds}

\medskip

In this section we prove theorem 1.2.A (see 2.2 below).

\medskip

{\it 2.1 A deformation of an almost complex structure.}

Recall that a finite set of vector fields
on a manifold $M$ satisfies {\it H\"ormander condition}
if these fields together with their iterated brackets
generate $TM$ at every point. A distribution 
(that is a field of tangent subspaces) on $M$ satisfies
H\"ormander condition if locally it can be generated
by such vector fields.

\proclaim Theorem 2.1.A. Let $(M,\Omega)$
be a closed symplectic manifold which admits an isotropic
distribution satisfying H\"ormander condition. Then
$\lambda_1(M,\Omega) = \infty.$

\medskip
\niN
{\bf Proof:} Let $J$ be an arbitrary compatible almost complex structure
on $M$. Denote by $g$ the corresponding Riemannian metric. 
Let $L$ be an isotropic distribution which satisfies
H\"ormander condition. Consider the following $g$-orthogonal
decomposition:
$$TM = L \oplus JL \oplus V.$$
Define a family of compatible almost complex structures $J_t$
on $M$ as follows: $J_t = t^{-1}J$ on $L$, $J_t = tJ$ on $JL$,
and $J_t = J$ on $V$. Notice that this deformation somewhat reminds
the Teichm\"uller deformation on surfaces. 

Let $g_t = t^{-1}g \oplus tg \oplus g$ be the Riemannian metric
associated with $J_t$. It suffices to show
that $\lambda_1(M,g_t)$ goes to infinity
when $t \to \infty$.

In order to prove this, choose vector fields
$X_1,...,X_N$ which are tangent to $L$ and satisfy H\"ormander
condition. Denote by $D_j$ the Lie derivative along $X_j$,
and by $D_j^*$ its $L^2$-adjoint. (The $L^2$- scalar product
on functions is associated to the canonical volume form on $M$).
Set $D = \Sigma_{j=1}^N D_j^* D_j$. The  H\"ormander theory [H]
implies
that $D$ has a discrete spectrum. Moreover, all eigenvalues
are non-negative, and the zero eigenspace consists of constant
functions only.

Let $H$ be the space of smooth 
functions on $M$ with the zero mean. The previous discussion
implies that $(Df,f) \geq C(f,f)$ for some positive constant
$C$ and all $f \in H$. Denote by $|\;|_t$ and $ \nabla_t$
the length on vectors  and the gradient in the sense of the metric $g_t$
respectively. Choose a constant $K >0$ such that $g(X_j,X_j) \leq K$
at every point for all $j=1,...,N$. Then

$$(Df,f) = \Sigma_{j=1}^N \int_M (D_jf)^2d{\rm Vol}
\leq \Sigma_{j=1}^N \int_M |\nabla_t f|^2_{t} |X_j|^2_t d{\rm Vol}$$
$$\leq NKt^{-1} 
\int_M |\nabla_t f|^2_t d{\rm Vol}.$$

Thus for some positive constant $C_1$ and for every $f \in H$
holds
$$
\int_M |\nabla_t f|^2_t d{\rm Vol} \geq C_1t \int_M f^2 d{\rm Vol}.$$
Therefore $\lambda_1(M,g_t) \to \infty$ when $t \to \infty$.
This completes the proof.
\QED

\medskip

\niN
{\bf Remark 2.1.B.} The idea of a collapse to a non-integrable
distribution
was extensively discussed in the literature on the first
eigenvalue, see e.g. [B-B],[Fu],[Ge].

\medskip

{\it 2.2 How to construct "generic" isotropic distributions?}

In general, existence of an isotropic distribution
imposes topological restrictions on $M$. However
this difficulty disappears after stabilization.

\proclaim Proposition 2.2.A. Let $(M,\Omega)$ be a closed
symplectic manifold. Then the stabilized manifold
$P = (T^4 \times M, \sigma \oplus \Omega)$ admits
a field of isotropic 2-planes which satisfies
H\"ormander condition.

\medskip
\niN
{\bf Proof of 1.2.A:} The theorem follows immediately
from 2.1.A and 2.2.A.
\QED

\medskip
\niN
{\bf Proof of 2.2.A:}
Clearly $P$ admits a 2-dimensional isotropic foliation.
A discussion in [G2] suggests that its generic perturbation
satisfies H\"ormander condition. We will use an explicit  argument.

1) Let $k$ be a positive integer. Consider $2k$ functions
$\phi_1(s),...,\phi_{2k}(s)$ defined as follows:
$\phi_{2j-1}(s) = \sin js$ and $\phi_{2j}(s) = \cos js$. 
Denote by $\phi^{(i)}$ the $i$-th derivative of a function $\phi$.
We claim that at every point $s$ the $2k \times 2k$ matrix 
$(\phi^{(i)}_j)$ is invertible. This is so since our functions
form a basis of solutions of a $2k$-th order differential
equation $$\Pi_{j=1}^k ({\partial^2 \over {\partial s^2}}
+ j^2)\phi = 0.$$

\medskip

2) Let $(p_1,q_1,p_2,q_2)$ be coordinates on $T^4$
such that the standard symplectic form is written as
$dp_1 \wedge dq_1 + dp_2 \wedge dq_2$.
Define an auxiliary set of vector fields $Z_j$ on $P$
as follows.
Set $Z_1 = 
{\partial \over {\partial p_1}}$ and
$Z_2 =
{\partial \over {\partial p_2}}$.
Let $Z_3,...,Z_{2k}$ be vector fields on $M$ which span $TM$
at every point (here $k$ is as large as needed).
Finally, set 
$$X_1 = 
{\partial \over {\partial q_1}} + \sin q_1 
{\partial \over {\partial p_2}},$$
and
$$X_2 =
{\partial \over {\partial q_2}}+
\Sigma_{j=1}^{2k} \phi_j(q_1)Z_j.$$

Clearly $X_1$ and $X_2$ generate an isotropic distribution
on $P$. Let us verify that it satisfies H\"ormander
condition. For that purpose denote $D = [X_1,.]$.
Then 
$$D^mX_2 = \Sigma_{j=1}^{2k} \phi_j^{(m)}(q_1)Z_j.$$

It follows from the step 1 of the proof that the span
of vector fields $$DX_2,...,D^{2k}X_2$$ at each point
contains vector fields $Z_1, ..., Z_{2k}$. Therefore
together with $X_1$ and $X_2$ these fields generate the whole $TP$.
This completes the proof.
\QED

\bigskip

{\bf 3. K-area, coupling and Hofer's norm}

\medskip

{\it 3.1 Two propositions}

Let $p: P \to S^2$ be a Hamiltonian symplectic fibration
with the fiber $(M,\Omega)$.
Theorem 1.9.A splits to 
the following two statements.

\medskip
\proclaim Proposition 3.1.A. The following identity
holds: $ \chi_+(P) = \size (P) = \epsilon (P)$.

\medskip

\proclaim Proposition 3.1.B. Suppose that the fibration $P$
is associated to an element $\gamma \in \pi_1(\Ham(M,\Omega))$.
Then
the following identities
hold: $\chi_+(P) = {1\over{||\gamma||_+}},$ and
$\chi(P) = {1\over{||\gamma||}}$. 

\medskip

The proofs are given in 3.2 and 3.3 respectively.

\vfill\eject

{\it 3.2 The coupling form}

Let $n$ be the complex dimension of $P$.
Recall [G-L-S],[McD-S]
that the {\it coupling form} $\delta^{\nu}$ of a
Hamiltonian symplectic
connection $\nu$ is the unique closed 2-form on $P$
whose restriction to each fiber coincides with the
symplectic form, which defines the
connection $\nu$ and such that the image of its top
power under the fiber integration map vanishes.
Here the fiber integration is a map, say $FI$ which
maps $2n$-forms on $P$ to 2-forms on $S^2$ as follows:

$$FI(\beta)(\xi,\eta) = \int_{P_{x}}i(\tilde \xi \wedge \tilde \eta)
\beta,$$
where $\xi,\eta$ are tangent vectors to $S^2$ at a point $x$,
and $\tilde \xi, \tilde \eta$ denote their (horisontal)
lifts to the points of the fiber $P_x$ over $x$.

It is proved in [G-L-S],[McD-S] that
for all $z \in P_x$ the coupling form can be written as
$\delta^{\nu}=
\Omega \oplus -\rho^{\nu} (z)$, where $TP = TP_x \oplus TS^2$
is the splitting associated to the connection $\nu$, and
$\rho^{\nu}(z)$ means the evaluation of the curvature form
at the point $z$ (remember that $\rho^{\nu}$ is a function-valued
form). Let us mention also, that this formula is cited in [P1,
proof of 1.6.B] with the opposite sign  at $\rho^{\nu}$.

\medskip
\niN
{\bf The sign convention.}
The Hamiltonian vector field $X$ generated by a Hamiltonian
function on $M$ is given by $dH = -i_X \Omega$. 

Denote by $a$ the positive generator of $H^2(S^2,\Z)$.
The proof of 3.1.A,B is based on the following simple fact.

\medskip

\proclaim Lemma 3.2.A.
Let $\omega$ be a symplectic form on $P$
whose restriction to fibers coincides with a 
a multiple of the symplectic form and
whose cohomology class equals to 
$p^*a+uc$ for some $u>0$.
Let $\nu$
be the symplectic connection defined by $\omega$-orthogonal
complements to the fibers.
Then $\omega - u\delta^{\nu} = p^*\tau$, where
$\tau$ is an area form on $S^2$ with the total area 1.

\medskip
\niN
{\bf Proof:}
Set $r=nu^{n-1}\int_M \Omega^{n-1}$, and define a 2-form
$\tau$ on $S^2$ by
$$\tau = {1\over r} FI(\omega^n).$$
Notice that the image of every  volume form, say $\beta$
on $P$ under
the fiber integration is always an area form
on $S^2$. Indeed for each positive frame $(\xi,\eta)$
on $S^2$ the form 
$i(\tilde \xi \wedge \tilde \eta)\beta$ does not vanishes
on the tangent space to the corresponding fiber. Thus $\tau$
is an area form.

It easily follows
from the definition of $FI$ above that $FI((\omega-p^*\tau)^n) = 0$.
In particular $(p^*a + uc -p^*[\tau])^n = 0$, thus the total
area of $\tau$ equals to 1. Moreover, the form
$u^{-1}(\omega - p^*\tau)$ satisfies the definition of the
coupling form of $\nu$. This completes the proof.
\QED

\vfill\eject
\niN
{\bf Proof of 3.1.A:}

1) Let $\omega$ be an arbitrary symplectic form on $P$
which extends symplectic structures on the fibers and defines
the positive orientation. Then $[\omega]= up^*a + c$ for some
$u>0$. Write $u^{-1}\omega = p^*\tau + u^{-1}\delta^{\nu}$,
where the area form $\tau$ and the connection $\nu$ are
defined  in 3.2.A. Decompose $\delta^{\nu}$ as $\Omega \oplus
-L^{\nu}\tau$. Then $u^{-1}\omega = \Omega \oplus (1-u^{-1}L^{\nu})
\tau$. Since this form is symplectic, we get that
$${1 \over {\max L^{\nu}}} > u^{-1}.$$
Note now that $s(\omega) = u^{-1}$. Therefore
$\chi_+(P) \geq \size (P)$.

\medskip

2) 
Take an arbitrary deformation of symplectic
forms satisfying (i) - (iv) of 1.8 , and take one of these forms
$\omega_u$. Set $\omega = u^{-1}\omega_u$. Then $\omega$
belongs to the class of symplectic forms which is used for the
definition of the size. Since $[\omega] = c + u^{-1}p^*a$ we
compute that ${\rm Vol} (P,\omega) = u^{-1} {\rm Vol} (M,\Omega)$,
thus $s(\omega) = u$. Taking $u$ arbitrarily close to $\epsilon (P)$
we conclude that $\size (P) \geq \epsilon (P)$.

\medskip

3) Let $\nu$ be a Hamiltonian symplectic connection on $P$.
Assume that for some area form $\tau$ of the total
area 1 holds ${1\over{\max L^{\nu}}} = r$.
Writing $\delta^{\nu} = \Omega \oplus -L^{\nu}\tau$
we see that
the deformation
$u\delta^{\nu} + \tau$ satisfies the conditions
(i)-(iv) of 1.8 provided $u < r$.
Hence $\epsilon(P) \geq r$, and taking $r$ arbitrarily
close to $\chi_+(P)$ we get the estimate
$\epsilon(P) \geq \chi_+(P)$. Together with the results
of the previous steps this completes the proof.
\QED

\medskip

{\it 3.3  Computing holonomy of a symplectic connection}

Here we prove 3.1.B.
The inequality
$\chi(P) \geq {1\over{||\gamma||}}$ was proved in [P1].
Exactly the same argument gives that
$\chi_+(P) \geq {1\over{||\gamma_+||}}.$ (One should be
careful with the sign conventions at this point!  )

So it remains to prove the opposite inequalities.
In what follows we will represent the 2-sphere $S^2$
as the unite square 
$$K = \{(x,y)|0\leq x\leq 1, 0\leq y \leq1\}$$
whose all boundary points are identified. We consider
Hamiltonian symplectic fibrations
over $K$ {\it trivialized over the boundary}, and Hamiltonian symplectic
connections with trivial parallel transport along every
segment of the boundary. Let $P$ be such a fibration,
and $\nu$ be such a connection on $P$. We write $M$
for the fiber over the boundary.

Fix an area form $\tau = dx \wedge dy$ on $K$.
Assume that the curvature of $\nu$ is given by $L(x,y,z)\tau$,
where
$(x,y) \in K$ and $z$ belongs to the fiber over $(x,y)$.

Consider a family of paths on $K$ given in $(x,y)$-coordinates
by
$\alpha(s) = \{(s,t)|t\in [0;1]\}$,
where $s\in [0;1]$. Note that $\{\alpha(s)\}$ represents
a {\it positive} generator of $\pi_1(LS^2)$ (here one
takes the boundary of $K$ as the base
point on $S^2$). Denote by $f_s :M \to M$ the holonomy of $\nu$
along $\alpha(s)$. Clearly, $\{f_s\}$ is a loop of Hamiltonian
diffeomorphisms of $M$. The needed inequalities follow
immediately from the next statement.

\medskip

\proclaim Lemma 3.3.A.
The following inequalities hold:
$$\length_+\{f_s\} \leq \max L$$ and
$$\length\{f_s\} \leq \max |L|.$$ 

\medskip
\niN
{\bf Proof:} Let $X,Y$ be the horizontal lifts
of the vector fields 
${{\partial } \over {\partial x}}$
and
${{\partial } \over {\partial y}}$
respectively. Let $X^s, Y^t$ be the flows
of $X$ and $Y$ respectively.
Set $h_{t,s} = Y^t X^s$ ({\it warning:} such
a map is defined on a domain which depends on
$t$ and $s$).

Note that $$(s,1,f^s(z)) = Y^1 X^s (0,0,z) = h_{1,s} (0,0,z). $$
Our goal is to calculate the vector field
$$v = 
{{\partial h_{1,s} } \over {\partial s}} (0,0,z).$$

Note that 
$${{\partial h_{t,s} } \over {\partial t}} = Y, \;
{{\partial h_{t,s} } \over {\partial s}} = Y^t_*X,$$
and therefore
$$
{{\partial } \over {\partial t}}
{{\partial h_{t,s} } \over {\partial s}} = 
{{\partial } \over {\partial s}}
{{\partial h_{t,s} } \over {\partial t}}-
[  
{{\partial h_{t,s} } \over {\partial t}},
{{\partial h_{t,s} } \over {\partial s}}]=
-[Y,Y^t_*X]=-Y^t_*[Y,X].$$

Notice that $[X,Y]$ is a vertical vector field
which equals to the Hamiltonian field
${\rm sgrad}L$, thus

$${{\partial } \over {\partial t}}
{{\partial h_{t,s} } \over {\partial s}} = 
{\rm sgrad} (L \circ Y^{-t}).$$

Taking into account that
$${{\partial h_{0,s} } \over {\partial s}} (0,0,z) = 
{{\partial } \over {\partial x}},$$ we get that
$$ v = 
{{\partial } \over {\partial x}}+
{\rm sgrad}\int_0^1 L \circ Y^{-t} dt\;,$$
and therefore the Hamiltonian function
generating the path $\{f_s\}$
is given by
$$F(s,z) = \int_0^1 L(s,1-t,Y^{-t}(s,1,z)) dt\;.$$
Clearly, $\max F \leq \max L$ and $\max |F| \leq \max |L|$.
This completes the proof of the lemma, and hence 
of the proposition 3.1.B.
\QED

\bigskip

{\bf 4. Estimating $\lambda_1$ on K\"ahler manifolds}

\medskip

In this section we prove Theorem 1.2.B. Below $a$ denotes
the positive generator of $H^2(\CP^N,\Z)$. We start with
two auxiliary results.

\medskip
{\it 4.1 Bourguignon-Li-Yau estimate} [B-L-Y]

Let $(M,\Omega,J,g)$ be a closed K\"ahler manifold
of complex dimension $n$, and let $\phi:M \to \CP^N$
be a holomorphic embedding. Then
$$\lambda_1(M,g) \leq 8\pi n
{{(\phi^*a \cup [\Omega]^{n-1},[M])}\over{([\Omega]^n,[M])}}.$$

\medskip

{\it 4.2 Demailly's theorem} [De, 8.5, p.58]

Let $(M,J)$ be a complex manifold of complex dimension $n$
with the canonical bundle $K$. There exists a constant $\alpha(n)$
which depends only on $n$ such that for every {\it ample} line
bundle $L \to M$, the bundle $$\alpha(n)(K+(n+2)L)$$
is {\it very ample}. In particular $M$ admits a holomorphic
embedding $\phi$ to a projective space with
$$\phi^*a = \alpha(n) (-c_1(TM) +(n+2)c_1(L)).$$

\medskip

{\it 4.3 Proof of Theorem 1.2.B:}

Let $J$ be a compatible complex structure on $(M,\Omega)$.
Since $[\Omega]$ is an integral class, there exists a Hermitian
holomorphic line bundle $L \to M$ whose Chern form equals to $\Omega$.
The Kodaira embedding theorem [Gr-Ha] implies that $L$ is ample.
Using 4.2 we get a holomorphic projective embedding $\phi$ of $M$
with $$\phi^*a = \alpha(n) (-c_1(TM) +(n+2)[\Omega]).$$
The required inequality 1.2.B follows now from 4.1.
\QED

\medskip

\niN
{\bf Acknowledgements.} My interest to 
the question about symplectic meaning
of the spectrum of the Laplacian arised after discussions
with Helmut Hofer during my stay at Bochum in 1991.
The relation between the coupling parameter
of a symplectic fibration  and the volume was suggested 
by Francois Lalonde and Dusa McDuff.
Paul Biran helped me with the theorem 1.2.B.
Dietmar Salamon explained me that 3.1.B is related to the 
Yang-Mills theory, and Thomas Davies gave me the reference
to Gravesen's paper [Gra].
And Yasha Eliashberg asked me a very stimulating question which pushed
me to extend some inequalities from [P1] to the identity 1.9.A.
I express my deep gratitude to all these people.

The work on this paper was completed during my visit to Lyon.
I thank Jean-Marie Morvan, Emmanuel Giroux and Albert Fathi
for their hospitality and useful conversations.

\bigskip

\end{document}